\documentclass[dvips,ejs]{imsart}

\newtheorem{proposition}{Proposition}[section]

\begin{document}

\begin{frontmatter}
\title{Inference in nonparametric current status models with covariates}
\runtitle{Current status models}

\begin{aug}
\author{\fnms{Odile} \snm{Pons}\ead[label=e1]{Odile.Pons@jouy.inra.fr}}
\address{INRA, Math\'ematiques, \\
78352 Jouy en Josas cedex, France\\
\printead{e1}} \runauthor{O. Pons}
\end{aug}

\begin{abstract}
In interval censored models with current status observations,
 the variables are indicators of the
presence of individuals on observation intervals and covariates.
When several individuals share the same
 observation interval, a simple procedure provides new estimators for
the distribution of the observation times and their intensity, in a
closed form. They are $n^{1/2}$-consistent for piece-wise constant
covariates. Estimators of the sample-sizes are deduced and
asymptotic $\chi^2$ tests for independence of the observations on
consecutive intervals and for independence between consecutive
classes for the observed individuals are proposed.
\end{abstract}


\tableofcontents
\end{frontmatter}

\section{Introduction}\label{introduction}
Statistical inference for sequential observations of individuals in
a large population differs according to the nature of the samples.
The observation of presence of individuals at specific locations is
often restricted to a sequence of time intervals. In
capture-recapture models, the size of finite and closed populations
has been estimated under the assumptions of the same parametric
model for the consecutive samples and time-dependent intensities for
the transitions of the populations between several states, with
individual covariates \cite{r1,r6,r7}.

The discrete observation sampling leads to cumulative observations
on fixed or random intervals, it is an interval censored model with
only current status observations. With individual observation times
for all the individuals, the monotonic nonparametric maximum
likelihood estimator of the time-dependent cumulative hazard
function relies on the greatest convex minorant algorithm, it weighs
the random observation times and converges at the rate $n^{1/3}$
(see \cite{r2,r3} and \cite{r4} in a model with constant
covariates). Here a nonparametric Markov model with piece-wise
constant covariate processes is considered as in \cite{r5} for
continuous observations, and the observations are current status
data with common observation intervals. A simple reparametrization
leads to easily calculated parametric estimators for the
distribution functions of the observation times and the population
sizes are estimated (section \ref{estimation}). The convergence
rates of the estimators in several nonparametric models is
$n^{1/2}$. In section \ref{dependent obs}, models with dependent
observations on consecutive time intervals are considered and new
estimators and tests for independence are proposed.

\section{Models with independent observations}\label{indep obs}

Consider a population of $L$ independent classes $C_1,\ldots, C_L$
of respective unknown sizes $\nu_l$, $l=1,\ldots, L$ and
$\nu=\nu_1+\ldots+\nu_L$. In each class, a sample of the population is
performed on a time interval $[0, \tau]$ with random sampling sizes $n_l$,
$l=1,\ldots, L$ and $n$.  Let $\tau_{l,1}<\ldots< \tau_{l,K_l}\leq \tau$
be the end-point observation intervals for class $C_l$ and
$(N_{li}(t))_{t\leq \tau}$ be the counting process of the observations
of individual $i$ of $C_l$ restricted to the intervals $I_{l,k}=
]\tau_{l,k-1}, \tau_{l,k}]$, $k=1, \ldots, K_l$ up to time $t$,
$$N_{li}(t)=\sum_{k=1}^{K_l} \delta_{li,k} 1\{I_{l,k}\cap[0, t]\neq
    \emptyset \}, \quad  \mbox{ with }
\delta_{li,k}= 1\{i\in C_l \mbox{ is observed on } I_{l,k} \},$$
with  $N_{li}(\tau)\leq K_l$, $\sum_{i=1}^{\nu_l}
1\{N_{li}(\tau)>0\}=n_l$. Only cumulated numbers $N_{li}(I_{l,k}$
are observed.

An individual $i$ of $C_l$ is supposed to be characterized by a
$p$-dimensional random covariate vector process $Z_{li}$ having
left-continuous sample-pathes with right-hand limits. The
individuals are sampled independently and for $l=1,\ldots, L$, the
processes $(N_{li}, Z_{li})$, $i=1,\ldots, n_l$, are mutually
independent and identically distributed. The distribution of
$N_{li}$ conditionally on  $Z_{li}$ is supposed to follow a Markov
model with independent increments, where the probability of
observing individuals only depends on their characteristics on the
observation interval
\begin{equation}
 \Pr(N_{li}(I_k) \vert (Z_{li}(s))_{s\leq \tau_{l,k}})=
\Pr(N_{li}(I_k) \vert Z_{li}(I_{l,k})), \label{Zk}
\end{equation}
only a countable set of values of the process $Z$ appears in the
whole sample-path of $N_{li}$.

The process $Z_{li}$ is sometimes restricted to a piece-wise constant
process with values $Z_{l,j}$ on a random sub-partition $I'_{li,j}=
[U_{li,j-1}, U_{li,j}[$, $j=1, \ldots, J$ of $(I_{l,k})_{l,k}$
\begin{equation}
Z_{li}(t)= \sum_{j=1}^{J} Z_{l,j} 1\{t\in I'_{li,j}\}.\label{Zj}
\end{equation}

The probability of observation of $i\in C_l$ on the partitions
$(I_{l,k})_k$ is a discrete process  defined according to the
assumption (\ref{Zk}) or (\ref{Zj}). Let $T_{li,k}$ be the unknown
first presence time of $i$ during the time interval $I_{l,k}$, and
we suppose that the model is defined by
\begin{eqnarray*}
p_{l,k}(Z_{li})&=& \Pr(\tau_{l,k-1}<  T_{li,k}\leq \tau_{l,k}\vert
Z_{li}) \\
&=& \sum_j  1\{t\in I'_{li,j}\subset I_{l,k}\}
\Pr(U_{li,j-1}<  T_{li,k}\leq U_{li,j}\vert Z_{li}(U_{li,j-1})) ,\\
P_{l}(Z_{l,j})&=& \Pr( Z_{li}(U_{li,j-1})= Z_{l,j}), \\
p_l&=& \Pr(N_{li}(\tau_{l,K_l})>0)= \int
\Pr(N_{li}(\tau_{l,K_l})>0\vert
Z_{li}(\tau_{l,K_l})) \, dP_{l}(Z_{li}) \\
&=& \sum_{j=1}^{J} \Pr(N_{li}(I'_{li,j}) >0\vert Z_{li}(U_{li,j-1})=
Z_{l,j}) \, P_{l}(Z_{l,j}) \\
&=&\sum_{k=1}^{K_l} \sum_{j=1}^{J} p_{l,k}(Z_{l,j})\,P_{l}(Z_{l,j}),\\
1-p_l&=& \Pr(N_{li}(\tau_{l,K_l})=0).
\end{eqnarray*}
However individuals $i$ with $N_{li}(\tau_{l,K_l})=0$ are not observed.
An underlying time-continuous model is defined by the intensities of
observation of the individuals.
The conditional intensity of observation of class $C_l$ is supposed to
    depend only on the current value of the covariate, for individual
    $i$ in $C_l$ and $t$ in $I_{l,k}$, it is defined by
$$\lambda_{l,k}(t, z)= \lim_{h\rightarrow 0}
    \frac{1}{h}\Pr(N_{li}(t+h)-N_{li}(t)>0\vert Z_{li}(t)=z)$$
More generally, the capture intensity for class $l$ is defined as
one of the intensity $\lambda_{l,k}$ by
\begin{eqnarray*}
 \lambda_{l}(t, Z_{li})&=& \lim_{h\rightarrow 0} \frac{1}{h}\Pr(
N_{li}(t+h)-N_{li}(t)>0\vert Z_{li}(t))\\
&=& \lim_{h\rightarrow 0}
\sum_{k=1}^{K_l} 1\{t\in I_{l,k}\} \sum_{j=1}^{J} 1\{t\in I'_{l,j}\}
 \lambda_{l,k}(t, Z_{l,j}) \mbox{ under (\ref{Zj})}.
\end{eqnarray*}
The variation of the cumulative intensities on each sub-interval are
denoted
\begin{eqnarray*}
\Delta\Lambda_{l,k}(t,Z_{li}) &=&\int_{I_{l,k}\cap [0,
    t]}\lambda_{l,k}(s,Z_{li}(s))\, ds \\
    &=&\sum_{j=1}^{J} 1\{I'_{li,j}\subset
I_{l,k}\} \int_{I'_{li,j}\cap [0, t]}\lambda_{l,k}(s,Z_{l,j})\, ds
\end{eqnarray*}
under (\ref{Zj}) and the cumulative intensities
from 0 is
$$\Lambda_{l}(t,Z_{li})=\sum_{k=1}^{K_l} 1\{t\in I_{l,k}\} \sum_{k'=1}^k
\Delta\Lambda_{l,k'}(t,z).$$ The unobserved apparition time
$T_{li,k}$ of $i$ in $C_l$ during the time interval $I_{l,k}$ has a
conditional distribution $\Pr\{T_{li,k}\leq t\vert
Z_{li}(\tau_{l,k})=Z_{l,j}\}=1-S_{l}(t, Z_{l,j})$, for a covariate
value $Z_{l,j}$. The probability of observation in $C_l$ is
continuously defined as
\begin{eqnarray*}
p_{l,k}(t,z)&=& \Pr(N_{li}(t)- N_{li}(\tau_{l,k-1})>0\vert
  Z_{li}(t)=z)= S_{l}(\tau_{l,k-1},z)- S_{l}(t,z)  \\
&=&  \exp\{- \Delta\Lambda_{l,k}(t,z)\}- \exp\{-
  \Delta\Lambda_{l,k}(\tau_{l,k-1},z)\}, \quad t\in I_{l,k},\\
p_l(t, Z_{li}) &=&\Pr(N_{li}(t)>0\vert  Z_{li})
 =1 - \exp\{- \Lambda_{l}(t, Z_{li}(t))\},
\end{eqnarray*}
$p_l(t, Z_{li})$ is the distribution function of observation for an
individual of $C_l$ before $t$ conditionally on the covariate. For $t$
in $I_{l,k}$, it is written
$p_l(t, Z_{li})= \sum_{k'<k} p_{l,k'}(Z_{li})  + p_{l,k}(t,Z_{li}).$\\

In a discrete nonparametric model, the hazard function of individual
$i$ in $C_l$ with covariate value $Z_{l,j}$ on an interval
$I'_{li,j}$ is written $\sum_k\lambda_{l,k}(t, Z_{l,j}) 1\{t\in
I_{l,k}\cap I'_{li,j}\}$.

The proportional hazards model is defined  by multiplicative intensities
$$\lambda_{l,k}(t, Z_{li})= \lambda_{l}(t) e^{\beta_{l,k}'Z_{li}(t)}=
\lambda_{l}(t) \sum_{j=1}^{J} e^{\beta_{l,k}'Z_{l,j}} 1\{t\in
I'_{li,j}\},$$ then
\begin{eqnarray*}
\Delta\Lambda_{l,k}(t, Z_{li}(t))&=& \sum_{j=1}^{J}
e^{\beta_{l,k}'Z_{l,j}} \int_{I_{l,k}\cap I'_{li,j}\cap[0,t]}
\lambda_{l}(s)\, ds\\
&=& \sum_{j=1}^{J} e^{\beta_{l,k}'Z_{l,j}} \Lambda_{l}(I_{l,k}\cap
I'_{li,j}\cap[0,t]).
\end{eqnarray*}

Let $S_{l}(t)= \exp\{-\int_0^t \lambda_{l}(s)\, ds\}$, for the
$\nu_l$ individuals, then the probability of being  unobserved is
$\Pr(T_{li}>\tau_{l,K_l})=  1 -p_l(\tau_{l,K_l})$, where $T_{li}$ the
first presence time of $i$,
\begin{eqnarray*}
 1 -p_l(\tau_{l,K_l}, Z_{li}) &=& \exp\{- \Lambda_{l}(\tau_{l,K_l},
 Z_{li})\}=  \prod_{k=1}^{K_l} \exp\{- \Delta\Lambda_{l,k}(\tau_{l,k},
 Z_{li})\}\\
&= & \prod_{k=1}^{K_l} \frac{S_{l}(\tau_{l,k},
 Z_{li})}{S_{l}(\tau_{l,k-1}, Z_{li})} \},\\
&=& \prod_{k=1}^{K_l}\prod_{j=1}^{J} 1\{I'_{li,j}\subset I_{l,k}\}
\{\frac{S_{l}(\tau_{l,k})}{S_{l}(\tau_{l,k-1})
 }\}^{\exp\{\beta_{l,k}'Z_{l,j}\}}
\end{eqnarray*}
and the conditional observation probability of $i$ on $I_{l,k}$ is
\begin{eqnarray*}
 p_{l,k}(Z_{li})&=& S_{l}(\tau_{l,k-1}, Z_{li}(\tau_{l,k}))-
 S_{l}(\tau_{l,k}, Z_{li}(\tau_{l,k})),\\
&=&\prod_{j=1}^{J} 1\{I'_{li,j}\subset I_{l,k}\}
\{\Delta S_{l}(I'_{li,j})\}^{\exp\{\beta_{l,k}'Z_{l,j}\}}.
\end{eqnarray*}

\section{Identifiability and estimation of the parameters}\label{estimation}
\subsection{Model without covariates}\label{scv}
Without covariates the parameters are only the probabilities
$p_{l,k}$ and $p_l(\tau_{l,K_l})$. Assuming that the observations on
the different intervals are independent, the model is multinomial
and the probabilities of independent observations on the $K_l+1$
intervals are written with the differences $\Delta_{l,k}=
\Delta\Lambda_{l,k}(I_{l,k}) >0, \,1\leq k\leq K_l$,
\begin{eqnarray}
&&1-p_l(\tau_{l,K_l})= \sum_{k\leq K_l} p_{l,k},\nonumber\\
&& \log (1-p_l(\tau_{l,k}))
= \sum_{k'\leq k} \{\log S_l(\tau_{l,k'-1}) - \log S_l(\tau_{l,k'})\}
=- \sum_{k'\leq k} \Delta_{l,k'},\label{pl,k}\\
&&\log p_{nl,k} = \log \{S_l(\tau_{l,k-1})- S_l(\tau_{l,k})\}
=\log \{1-\exp(-\sum_{k'\leq k} \Delta_{l,k'})\}.
\nonumber
\end{eqnarray}
 The log-likelihood for class $C_l$ is
$$l_n(l)= \sum_{i=1}^{n_l}  [\sum_{k\leq K_l}  \{\delta_{li,k} \log p_{l,k}
+ (1- \delta_{li,k}) \log  (1-p_{l,k})\} ]$$ under (\ref{Zj})  and
the MLE of the parameters $p_{l,k}$ and the function $S_{l}$ are
\begin{eqnarray*}
\widehat{p}_{nl,k}&=& n_l^{-1} \sum_{i=1}^{n_l} \delta_{li,k},\quad
\widehat{p}_{nl}(\tau_{l,K_l})= 1- n_l^{-1} \sum_{i=1}^{n_l}
\sum_{k=1}^{K_l}\delta_{li,k},\\
\widehat{S}_{nl}(\tau_{l,k})&=& \widehat{S}_{nl}(\tau_{l,k-1})-
\widehat{p}_{nl,k} = 1-n_l^{-1} \sum_{i=1}^{n_l} \sum_{k'=1}^{k}
\delta_{li,k'}.
\end{eqnarray*}
The estimator $\widehat{S}_{nl}$ is decreasing with weights at the
sampling times $\tau_{l,k}$. From (\ref{pl,k}), the differences
$\Delta_{l,k}$ satisfy
$$\Delta_{l,k}= \log \frac{1-\sum_{k'<k} p_{l,k}}{1-\sum_{k'\leq k}
p_{l,k}} >0,$$ their estimators are deduced from the
$\widehat{p}_{nl,k}$'s and the cumulative hazard function for $C_l$ is
estimated by
\begin{equation}
 \widehat{\Lambda}_{nl}(t)=\sum_{k=1}^K 1\{\tau_{l,k-1}<t\leq
 \tau_{l,k}\}\log \frac{1-\sum_{k'<k} \widehat{p}_{nl,k}}{1-\sum_{k'\leq k}
 \widehat{p}_{nl,k}}.\label{Lambdan}
\end{equation}

Let $p_{0l,k}$, $S_{0l}$ and $\Lambda_{0l}$ be the actual values of the
model parameters, then
\begin{proposition}\label{cv1}
The estimators $\widehat{p}_{nl,k}$, $\widehat{\Lambda}_{nl,k}$ and
$\widehat{S}_{nl}$ are a.s. consistent as $n\rightarrow\infty$,
$n_l^{1/2}(\widehat{p}_{nl,k}- p_{0l,k})_{k}$ converge to
centered Gaussian variable with covariances
$n_l^{-1}p_{0l,k}(1-p_{0l,k})$ and zero otherwise,
and the processes $n_l^{1/2}(\widehat{S}_{nl}-S_{0l})$ and
$n_l^{1/2}(\widehat{\Lambda}_{nl}-\Lambda_{0l})$ converge to
centered Gaussian process with independent increments and
variances
\begin{eqnarray*}
n_lE(\widehat{S}_{nl}-S_{0l})^2(\tau_{l,k}) &=& \sum_{k'<k}
p_{0l,k'}(1-p_{0l,k'}),\\
n_lE(\widehat{\Lambda}_{nl}-\Lambda_{0l})^2(\tau_{l,k}) &=&
\sum_{k'<k} p_{0l,k'}(1-p_{0l,k'})
\left(\frac{p_{0l,k}}{\Pr(T_{li}>\tau_{l,k-1})
\Pr(T_{li}>\tau_{l,k})} \right)^2 \\
&&+p_{0l,k}(1-p_{0l,k})\left(\frac{1}{\Pr(T_{li}>\tau_{l,k})}
\right)^2.
\end{eqnarray*}
\end{proposition}

\subsection{Models with covariates}\label{covariatemodel}
The  parameters of the model are the probabilities $p_l$ and
 $p_{l,k}=p_l(I_{l,k})$, or the functions $p_l(z)$ and
 $p_{l,k}(z)=p_l(I_{l,k},z)$ in regression model. The probabilities $p_l$ are
 expressions of the  $p_{l,k}$'s and of the distribution of the
 covariates, their estimators satisfy
\begin{eqnarray}
\widehat{p}_{nl,k}&=& \sum_{j=1}^J \widehat{p}_{nl,k}(Z_{l,j})\,
\widehat{p}_{nl}(Z_{l,j}),\label{p2}\\
\widehat{p}_l&=& \sum_{k=1}^{K_l}\sum_{j=1}^J
\widehat{p}_{nl,k}(Z_{l,j}) \,\widehat{p}_{nl}(Z_{l,j})\nonumber
\end{eqnarray}
but the distributions $p_l$ are not directly estimable since all
the individuals are not observed. Only the probabilities
$\Pr(Z_{li}\leq z\vert \delta_{li,k}=1)$ are directly estimable as the
proportion of the individuals  observed in $I_{lk}$ such that
$Z_{li}\leq z$.  Then $P_{l}(z)$ is deduced from the equation
\begin{equation}
P_{l}(z) = \frac{\sum_{j=1}^J \Pr(Z_{li}\leq z\vert \delta_{li,k}=1)
\Pr(\delta_{li,k}=1)}{\sum_{j=1}^J \Pr(\delta_{li,k}=1 \vert
  Z_{li}\leq z)},\, \forall i= 1,\, \ldots, n \label{pobs}
\end{equation}
which is easily estimated with the empirical probabilities.\\

The estimable parameters are always the values of the functions $S_{l}$ and
 $\Lambda_{l}$ at the observation times $\tau_{l,k}$ and model
 parameters when it is appropriate.
Conditionally on the covariates, the log-likelihood for class $C_l$ is
\begin{eqnarray*}
l_n(l)&=& \sum_{i=1}^{n_l}\sum_{k\leq K_l}\{\delta_{li,k} \log p_{l,k}(Z_{li})
+ (1- \delta_{li,k}) \log  (1-p_{l,k}(Z_{li}))\}\\
&=& \sum_{i=1}^{n_l}  \sum_{k\leq K_l} \sum_{j=1}^{J} 1\{I'_{li,j}\subset
  I_{l,k}\} \{\delta_{li,k} \log p_{l,k}(Z_{l,j})\\
&&+ (1- \delta_{li,k}) \log  (1-p_{l,k}(Z_{l,j}))\}.
\end{eqnarray*}
The MLEs are identical to the previous estimators if the covariates are
on the intervals  $I_{l,k}$ and $p_{l,k}(Z_{li})\equiv p_{l,k}$. If $J$ is
finite, and the variations of the processes $Z_{l,i}$ are observed
though those of $N_{l,i}$ are only observed on $I_{l,k}$, $i=1,\ldots,
n$, they are modified
\begin{eqnarray*}
\widehat{p}_{nl,k}(Z_{l,j})&=& n_l^{-1} \sum_{i=1}^{n_l}
\delta_{li,k} 1\{I'_{li,j}\subset I_{l,k}\},\\
\widehat{S}_{nl}(\tau_{l,k}, Z_{l,j})&=& 1-n_l^{-1} \sum_{i=1}^{n_l}
\sum_{k'=1}^{k}  \delta_{li,k'} 1\{I'_{li,j}\subset I_{l,k}\},\\
\widehat{S}_{nl}(\tau_{l,k},z)&=& 1-n_l^{-1} \sum_{i=1}^{n_l} \sum_{k'=1}^{k}
\delta_{li,k'} \sum_{j=1}^{J} 1\{Z_{l,j}=z\}1\{I'_{li,j}\subset I_{l,k}\},\\
\widehat{\Lambda}_{nl}(t,z)&=& \sum_{k=1}^K \sum_{j=1}^{J} 1\{t\in
I'_{li,j}\subset I_{l,k}\} 1\{Z_{l,j}=z\}
\log \frac{1-\sum_{k'<k} \widehat{p}_{nl,k}(z)}{1-\sum_{k'\leq k}
 \widehat{p}_{nl,k}(z)}.
\end{eqnarray*}
\, \\
With continuous covariate and under (\ref{Zk}), kernel estimators of
the functions conditionally on $z$ are defined with a kernel $K$, a
bandwidth $h$ and $K_h(x)=h^{-1} K(h^{-1}x)$, by smoothing these
estimators or the previous ones
\begin{eqnarray*}
\widehat{p}_{nl,k}(z)&=& \frac{\sum_{i=1}^{n_l} K_h(z-Z_{l,i}(\tau_{l,k}))
  \delta_{li,k}}{\sum_{i=1}^{n_l} K_h(z-Z_{l,i}(\tau_{l,k}))},\\
\widehat{S}_{nl}(\tau_{l,k}, z)&=&1-\sum_{k'=1}^{k}\widehat{p}_{nl,k'}(z),\\
\widehat{\Lambda}_{nl}(t,z)&=&\sum_{k=1}^K  \frac{\sum_{i=1}^{n_l}
  K_h(z-Z_{l,i}(\tau_{l,k}))\delta_{li,k}}{\sum_{i=1}^{n_l}
  K_h(z-Z_{l,i}(\tau_{l,k}))}\\
  &&\qquad \qquad \times
\sum_{j=1}^{J} 1\{t\in I_{l,k}\} \log \frac{1-\sum_{k'<k}
  \widehat{p}_{nl,k}(z)}{1-\sum_{k'\leq k} \widehat{p}_{nl,k}(z)}
\end{eqnarray*}
and they converge at the usual rate of the kernel estimators if the
bandwidth tends to zero at the optimal rate $n^{-\frac{s}{d+4s}}$,
for a $p$-dimensional covariate having a density with a s-order derivative. \\

For estimation in the proportional hazards model with constant
covariates $Z_{li,k}$ on $I_{i,k}$, let $\omega_{li,k}= \exp\{\beta'_{l,k}
Z_{li,k}\}$, $\Omega_l=\{\omega_{li,k}\}_{i\leq n k\leq K_l}$,
\begin{eqnarray}
&&\log\Delta S_{l}(I_{l,k})= \log S_{l}(\tau_{l,k-1})+ \log \{1-
 \frac{S_{l}(\tau_{l,k})}{S_{l}(\tau_{l,k-1})}\}\nonumber\\
 &&\qquad = \sum_{k'\leq k}\Delta_{l,k'}-\log  (1-e^{-\Delta_{l,k}}),\nonumber\\
&& \log (1-p_l(\tau_{l,K_l}, \omega_{li,k}))
 =-\sum_{k=1}^{K_l} \omega_{li,k} \Delta_{l,k}, \label{pDelta}\\
&&\log p_{l,k}(Z_{li}) = \omega_{li,k} \log\Delta S_{l}(I_{l,k}) =-
\omega_{li,k} \{\sum_{k'<\leq k}\Delta_{l,k'}-\log
(1-e^{-\Delta_{l,k}})\}\nonumber.
\end{eqnarray}
Denote $\mu_{l,k}= \log\Delta S_{l}(I_{l,k})= \log p_{l}(I_{l,k})$,
then the estimator of $p_{l,k}(Z_{li,k})=
\exp\{\omega_{li,k}\mu_{l,k}\}$ of proposition \ref{cv1} has to be
restricted to the individuals with the same covariate value as
$Z_{li,k}$.
\begin{proposition}
If $\Omega_l$ is a finite set $\{\omega_{l,j}\}_{j=1,\ldots,J}$, then
$$\omega_{l,j}= \log\frac{p_{l}(I_{l,k}, Z_{li,k})}{p_{l}(I_{l,k})},$$
and estimators are defined by
\begin{eqnarray*}
\widehat{p}_{nl}(I_{l,k}, Z_{l,j})&=&  \frac{\sum_{i\leq n_l}
  1\{\omega_{li,k}= \omega_{l,j}\}\delta_{li,k}} {\sum_{i\leq n_l}
  1\{\omega_{li,k}=\omega_{l,j}\}},\\
\widehat{\mu}_{nl,k}&=& \log\widehat{p}_{nl,k}= \log\{n_l^{-1}
\sum_{i=1}^{n_l} \delta_{li,k}\},\\
\widehat{\omega}_{nl,j}&=& \log \frac{n_l(\sum_{i\leq n_l}1\{\omega_{li,k}=
  \omega_{l,j}\}\delta_{li,k})}{(\sum_{i\leq n_l}\delta_{li,k})
(\sum_i1\{\omega_{li,k}=\omega_{l,j}\})},\\
\widehat{S}_{nl}(\tau_{l,k}, Z_{l,j})&=& 1- \frac{\sum_{i\leq
  n_l}\sum_{k'=1}^k   1\{\omega_{li,k}= \omega_{l,j}\}\delta_{li,k}}
  {\sum_{i\leq n_l} 1\{\omega_{li,k}=\omega_{l,j}\}}.
\end{eqnarray*}
\end{proposition}
An estimator of $\Lambda_{l}(\tau_{l,k}, Z_{l,j})$ is deduced from
the $\widehat{p}_{nl}(I_{l,k}, Z_{l,j})$'s and (\ref{pl,k}) as
previously, $$ \widehat{\Lambda}_{nl}(t, Z_{l,j})=\sum_{k=1}^K
1\{\tau_{l,k-1}<t\leq
 \tau_{l,k}\}\log \frac{1-\sum_{k'<k} \widehat{p}_{nl}(I_{l,k}, Z_{l,j})}
 {1-\sum_{k'\leq k} \widehat{p}_{nl}(I_{l,k}, Z_{l,j})}$$
 and the results of Proposition \ref{cv1} extend to these estimators.

Let $p_{0l,k}$, $S_{0l}$ and $\Lambda_{0l}$ be the actual values of
the model parameters, then
\begin{proposition}
The estimators $\widehat{p}_{nl,k}$, $\widehat{\Lambda}_{nl,k}$ and
$\widehat{S}_{nl}$ are a.s. consistent as $n\rightarrow\infty$,
$n_l^{1/2}(\widehat{p}_{nl,k}- p_{0l,k})_{k}$ converge to centered
Gaussian variable with covariances $n_l^{-1}p_{0l,k}(1-p_{0l,k})$
and zero otherwise, and the processes
$n_l^{1/2}(\widehat{S}_{nl}-S_{0l})$ and
$n_l^{1/2}(\widehat{\Lambda}_{nl}-\Lambda_{0l})$ converge to
centered Gaussian process with independent increments and variances
\begin{eqnarray*}
n_lE(\widehat{S}_{nl}-S_{0l})^2(\tau_{l,k}) &=&  \sum_{k'<k}
p_{0l,k'}(1-p_{0l,k'}),\\
n_lE(\widehat{\Lambda}_{nl}-\Lambda_{0l})^2(\tau_{l,k}) &=&
\sum_{k'<k} p_{0l,k'}(1-p_{0l,k'})
\left(\frac{p_{0l,k}}{\Pr(T_{li}>\tau_{l,k-1})
\Pr(T_{li}>\tau_{l,k})} \right)^2 \\
&&+ p_{0l,k}(1-p_{0l,k})\left(\frac{1}{\Pr(T_{li}>\tau_{l,k})}
\right)^2.
\end{eqnarray*}
\end{proposition}

The proportional hazards model without finite $\Omega_{l}$  is still
parametric but maximum likelihood estimators are not written  in closed
form. Denoting $\Delta_{li,j}=  \Lambda_{l}(U_{li,j})-
\Lambda_{l}(U_{li,j-1})$, the probabilities are now
\begin{eqnarray*}
&&\log (1-p_l(\tau_{l,K_l}, \beta_{k,l}, Z_{l,j})) \\
&&\quad=\sum_{k=1}^{K_l}\sum_{j=1}^{J} 1\{I'_{li,j} \subset
I_{l,k}\}
\exp\{\beta_{l,k}'Z_{l,j}\} \{\log S_{l}(U_{li,j})-\log S_{l}(U_{li,j-1})\}\\
&& \quad =-\sum_{k=1}^{K_l}\sum_{j=1}^{J} 1\{I'_{li,j}\subset I_{l,k}\}
\exp\{\beta_{l,k}'Z_{l,j}\} \Delta_{li,j}, \\
&&\log p_{l,k}(Z_{li}) = \sum_{j=1}^{J} 1\{I'_{li,j}\subset I_{l,k}\}
\exp\{\beta_{l,k}'Z_{l,j}\} \log\Delta S_{l}(I'_{li,j})\\
&& \quad =\sum_{j=1}^{J} 1\{I'_{li,j}\subset I_{l,k}\}
\exp\{\beta_{l,k}'Z_{l,j}\} [\log S_{l}(U_{li,j-1}) + \log \{1-
 \frac{S_{l}(U_{li,j})}{S_{l}(U_{li,j-1})}\}]\\
&& \quad =- \sum_{j=1}^{J} 1\{I'_{li,j}\subset
I_{l,k}\}\exp\{\beta_{l,k}'Z_{l,j}\} [\sum_{j'<
j}\Delta_{li,j'}\\&&\qquad\qquad\qquad \qquad\qquad+ \log
\{1-\exp(-\exp\{\beta_{l,k}'Z_{l,j}\} \Delta_{li,j})\}].
\end{eqnarray*}
When covariate only depend on the observation intervals, the parameters
are all identifiable by maximization of the likelihood, as it is the
case with continuously observed individuals. The parameters are not
identifiable when the covariates vary individually.

\subsection{Estimation of the sample size}\label{estim size}
The unknown population size $\nu$ has to be estimated. For a
population of $L$ observed classes $C_1,\ldots, C_L$ of respective
sizes $\nu_l$, estimators of the catching or observation
probabilities $p_{l,k}$ would be $n_{l,k} \nu_{l}^{-1}$ if $\nu_{l}$
was known, $k=1,\ldots, K_l$. By inverting this expression after an
estimator $\widehat{p}_{nl}$ has been defined, the sizes are usually
estimated by
$$ \widehat{\nu}_{nl} = \frac{n_l}{\widehat{p}_{nl}},\, l=1,\ldots, L, \qquad
\widehat{\nu}_n=\sum_{l=1}^L \widehat{\nu}_{nl} =\sum_{l=1}^L
\frac{n_l}{\widehat{p}_{nl}}.$$
With consecutive intervals under the same conditions and with varying
catching or observation probabilities $p_{l,k}$, define a moving
average estimator of $p_{l,k}$ and mean estimators of classes and
population sizes for $k>a\geq 1$ by
$$ \widehat{p}_{nl, k}= \frac{\sum_{k'=k-a}^{k+a} \widehat{p}_{nl,
k'}}{2a},\quad \widehat{\nu}_{nl}=
\sum_{k>a}\frac{n_{l,k}}{\widehat{p}_{nl, k}},\quad \widehat{\nu}_n=
\sum_{l=1}^L \widehat{\nu}_{nl}.$$ The same method applies for
covariate dependent probabilities, using the estimators of section
\ref{covariatemodel} and (\ref{p2})-(\ref{pobs}).

\section{Models with dependent observations on consecutive intervals}
\label{dependent obs}
\subsection{Nonparametric models}\label{npmodel}
When the probability of observing individuals in $I_{l,k}$ depends
on their observation in  $I_{l,k-1}$, several nonparametric models
may be considered. Let
\begin{eqnarray*}
 \pi_{l,k}&=& \Pr\{\tau_{l,k-1}<T_{li}\leq  \tau_{l,k+1}\vert
\tau_{l,k-1}<T_{li}\leq  \tau_{l,k}\},\\
\pi_{l,k}(Z_{li})&=& \Pr\{\tau_{l,k-1}<T_{li}\leq  \tau_{l,k+1}\vert
\tau_{l,k-1}<T_{li}\leq  \tau_{l,k}, Z_{li}\},
\end{eqnarray*}
then
$$p_{l,k, k+1}= \Pr\{\tau_{l,k-1}<T_{li}\leq  \tau_{l,k+1}\} =
\pi_{l,k} p_{l,k}$$
and conditionally on $Z_{li}$, $p_{l,k, k+1}(Z_{li})=
\pi_{l,k}(Z_{li}) p_{l,k}(Z_{li})$.
The estimators are now defined for joint intervals,
\begin{eqnarray*}
\widehat{\pi}_{nl,k}&=& \frac{\sum_{i=1}^{n_l}
\delta_{li,k}\delta_{li,k+1}} {\sum_{i=1}^{n_l}\delta_{li,k}},\\
\widehat{p}_{nl,k,k+1}&=&n_l^{-1}\sum_{i=1}^{n_l}\delta_{li,k}
\delta_{li,k+1},\\
\widehat{p}_{nl,k, k+1}(Z_{l,j})&=& n_l^{-1} \sum_{i=1}^{n_l} \delta_{li,k}
\delta_{li,k+1} 1\{I'_{li,j} \subset I_{l,k}\cup I_{l,k+1}\},\\
\widehat{\pi}_{nl,k}(Z_{l,j})&=& \frac{\sum_{i=1}^{n_l}
\delta_{li,k}\delta_{li,k+1}1\{I'_{li,j} \subset I_{l,k}\cup I_{l,k+1}\}}
{\sum_{i=1}^{n_l}\delta_{li,k}1\{I'_{li,j} \subset I_{l,k}\}}.
\end{eqnarray*}
All the other models and estimators of section \ref{npmodel} are
generalized by the same method. In the model without covariates, a
test for the hypothesis $H_0$ of independence between intervals
$I_{l,k}$ and $I_{l,k+1}$ is a test for $p_{l,k,
k+1}=p_{l,k}p_{l,k+1}$ or $\pi_{l,k}= p_{l,k+1}$.
\begin{proposition}
Under  $H_0$, the statistic
$$Z_l=\sum_{k=1}^{K_l-1}\frac{(\widehat{p}_{nl,k} \widehat{p}_{nl,k+1}-
  \widehat{p}_{nl,k, k+1})^2}{ \widehat{p}_{nl,k}\widehat{p}_{nl,k+1}}$$
converges to a $\chi^2_{(K_l-2)^2}$ as $n_l\rightarrow\infty$.
\end{proposition}
{\it Proof}. Let $N_{l,k} =\sum_{i=1}^{n_l} N_{li}(I_{l,k})$, $N_{l,k,k+1}
=\sum_{i=1}^{n_l} N_{li}(I_{l,k}\cup I_{l,k+1})$ and
$$Z_l=\sum_{k=1}^{K_l-1}\frac{( N_{l,k,k+1}- n_l^{-1}N_{l,k}
  N_{l,k+1})^2}{N_{l,k} N_{l,k+1}}$$
is the test statistic for independent marginals in a two-dimensional array.

\subsection{Markov models}\label{markovmodel}
As the individual classes change during the observation period, a
second class index may be incorporated in the model to take into
account the evolution. Let $C_{i,T_i}$ denote the class at $T_i$ for
some observation time $T_i$ of individual $i$,
\begin{eqnarray*}
\eta_{ll',i} &=&1\{C_{i,T_i}=C_l,  C_{i,T^-_i}=C_{l'}\},\\
 p_{l\vert l',k}&=&p_{l\vert l'}(I_{l,k})= \Pr\{T_i\in I_{l,k},
 C_{i,T_i}=C_l\vert  C_{i,T^-_i}=C_{l'}\},\\
S_{l\vert l',k}&=&\Pr\{T_i\in I_{l,k},T_i\geq t, C_{i,T_i}=C_l\vert
 C_{i,T^-_i}=C_{l'}\},\\
\Lambda_{l\vert l',k}&=&h^{-1}\lim_{h\rightarrow 0}\Pr\{T_i\in [t, t+h),
 C_{i,T_i}=C_l\vert T_i\geq t, C_{i,T^-_i}=C_{l'}\},
\end{eqnarray*}
The likelihood is proportional to
\begin{eqnarray*}
\prod_{l=1}^{L} \prod_{k=1}^{K_l}  \prod_{i=1}^n \prod_{l'=1}^L
 \{p_{l\vert l',k}^{\delta_{li,k}}
(1-p_{l\vert l',k})^{1-\delta_{li,k}}\}^{\eta_{ll',i}}
\end{eqnarray*}
and the estimators become
\begin{eqnarray*}
\widehat{p}_{nl\vert l',k}&=& \frac{\sum_{i=1}^n
  \delta_{li,k}\eta_{ll',i}}{\sum_{i=1}^{n_l} \eta_{ll',i}},\\
\widehat{S}_{nl\vert l'}(\tau_{l,k})&=& 1- \frac{\sum_{i=1}^{n_l}
\sum_{k'=1}^{k} \delta_{li,k'}\eta_{ll',i}}{\sum_{i=1}^{n_l} \eta_{ll',i}},\\
 \widehat{\Lambda}_{nl\vert l'}(t)&=& \sum_{k=1}^K \sum_{j=1}^{J}
1\{\tau_{l,k-1}<t\leq \tau_{l,k}\}
\log \frac{1-\sum_{k'<k} \widehat{p}_{nl\vert l',k}}{1-\sum_{k'\leq k}
  \widehat{p}_{nl\vert l',k}}.
\end{eqnarray*}
The extension to models and estimators with covariates follows
easily from section \ref{covariatemodel}. A test for the hypothesis
$H_0$ of independence between observation and the variation between
classes is a test for $p_{l\vert l',k}=p_{l,k}\Pr\{C_{i,T_i}=C_l|
C_{i,T^-_i}=C_{l'}\}$ for every $l, l'=1,\ldots, L$ and $k=1,\ldots,
K_l$.

Let $q_{ll'}=\Pr\{C_{i,T_i}=C_l, C_{i,T^-_i}=C_{l'}\}$, then the
estimators
$$\widehat{q}_{nl l'}=\frac{\sum_{i=1}^{n_l}
\eta_{ll',i}}{n_l},\quad \widehat{p}_{nl l',k}= \widehat{p}_{nl|
l',k}\, \widehat{q}_{nl l'},$$ provide a test statistic.
\begin{proposition}
Under  $H_0$, the statistic
$$X_l=\sum_{k=1}^{K_l} \sum_{l=1}^{L}
\frac{(\widehat{p}_{nll',k}-\widehat{p}_{nl,k} \,
\widehat{q}_{nll'})^2}{\widehat{p}_{nl,k} \,\widehat{q}_{nll'}}$$
converges to a $\chi^2_{(K_l-1)(L-1)}$ as $n_l\rightarrow\infty$.
\end{proposition}

\end{document}